\documentclass[12pt]{article}
\usepackage{amsfonts,amsmath,amssymb,setspace,mathtools,bm}
\usepackage{graphicx,float,cleveref}
\usepackage{tikz}
\usepackage{amsthm}
\usepackage{fullpage}
\usepackage{ifthen}             
\usetikzlibrary{calc}           

\def\rotateclockwise#1{
  \newdimen\xrw
  \pgfextractx{\xrw}{#1}
  \newdimen\yrw
  \pgfextracty{\yrw}{#1}
  \pgfpoint{\yrw}{-\xrw}
}

\def\rotatecounterclockwise#1{
  \newdimen\xrcw
  \pgfextractx{\xrcw}{#1}
  \newdimen\yrcw
  \pgfextracty{\yrcw}{#1}
  \pgfpoint{-\yrcw}{\xrcw}
}

\def\outsidespacerpgfclockwise#1#2#3{
  \pgfpointscale{#3}{
    \rotateclockwise{
      \pgfpointnormalised{
        \pgfpointdiff{#1}{#2}}}}
}

\def\outsidespacerpgfcounterclockwise#1#2#3{
  \pgfpointscale{#3}{
    \rotatecounterclockwise{
      \pgfpointnormalised{
        \pgfpointdiff{#1}{#2}}}}
}

\def\outsidepgfclockwise#1#2#3{
  \pgfpointadd{#2}{\outsidespacerpgfclockwise{#1}{#2}{#3}}
}

\def\outsidepgfcounterclockwise#1#2#3{
  \pgfpointadd{#2}{\outsidespacerpgfcounterclockwise{#1}{#2}{#3}}
}

\def\outside#1#2#3{
  ($ (#2) ! #3 ! -90 : (#1) $)
}

\def\cornerpgf#1#2#3#4{
  \pgfextra{
    \pgfmathanglebetweenpoints{#2}{\outsidepgfcounterclockwise{#1}{#2}{#4}}
    \let\anglea\pgfmathresult
    \let\startangle\pgfmathresult

    \pgfmathanglebetweenpoints{#2}{\outsidepgfclockwise{#3}{#2}{#4}}
    \pgfmathparse{\pgfmathresult - \anglea}
    \pgfmathroundto{\pgfmathresult}
    \let\arcangle\pgfmathresult
    \ifthenelse{180=\arcangle \or 180<\arcangle}{
      \pgfmathparse{-360 + \arcangle}}{
      \pgfmathparse{\arcangle}}
    \let\deltaangle\pgfmathresult

    \newdimen\x
    \pgfextractx{\x}{\outsidepgfcounterclockwise{#1}{#2}{#4}}
    \newdimen\y
    \pgfextracty{\y}{\outsidepgfcounterclockwise{#1}{#2}{#4}}
  }
  -- (\x,\y) arc [start angle=\startangle, delta angle=\deltaangle, radius=#4]
}

\def\corner#1#2#3#4{
  \cornerpgf{\pgfpointanchor{#1}{center}}{\pgfpointanchor{#2}{center}}{\pgfpointanchor{#3}{center}}{#4}
}

\def\hedgeiii#1#2#3#4{
  \outside{#1}{#2}{#4} \corner{#1}{#2}{#3}{#4} \corner{#2}{#3}{#1}{#4} \corner{#3}{#1}{#2}{#4} -- cycle
}

\def\hedgem#1#2#3#4{
  
  \outside{#1}{#2}{#4}
  \pgfextra{
    \def\hgnodea{#1}
    \def\hgnodeb{#2}
  }
  foreach \c in {#3} {
    \corner{\hgnodea}{\hgnodeb}{\c}{#4}
    \pgfextra{
      \global\let\hgnodea\hgnodeb
      \global\let\hgnodeb\c
    }
  }
  \corner{\hgnodea}{\hgnodeb}{#1}{#4}
  \corner{\hgnodeb}{#1}{#2}{#4}
  -- cycle
}

\newtheorem{Theorem}{Theorem}
\newtheorem{Lemma}[Theorem]{Lemma}

\newtheorem{Question}{Question}

\newenvironment{Proof}{\begin{trivlist} \item[] {\bf Proof.}}{\hfill $\Box$\end{trivlist}}

\renewcommand{\geq}{\geqslant}
\renewcommand{\leq}{\leqslant}

\tikzstyle{vertex} = [fill,shape=circle,node distance=80pt]
\tikzstyle{edge} = [fill,opacity=.5,fill opacity=.5,line cap=round, line join=round, line width=50pt]
\tikzstyle{elabel} =  [fill,shape=circle,node distance=30pt]
\pgfdeclarelayer{background}
\pgfsetlayers{background,main}
\usetikzlibrary{graphs,graphs.standard}

\title{Existential Closure in Uniform Hypergraphs}
\author{Andrea C. Burgess\thanks{Department of Mathematics and Statistics, University of New Brunswick, Saint John, NB, E2L 4L5, Canada, andrea.burgess@unb.ca}
       \and Robert D. Luther\thanks{Department of Mathematics and Statistics, Memorial University of Newfoundland, St. John's, NL, A1C 5S7, Canada, rdl863@mun.ca (corresponding author)}
       \and David A. Pike\thanks{Department of Mathematics and Statistics, Memorial University of Newfoundland, St. John's, NL, A1C 5S7, Canada, dapike@mun.ca}}

\begin{document}
\begin{doublespace}
\date{\today}

\maketitle

\begin{abstract}
For a positive integer $n$, a graph with at least $n$ vertices is $n$-existentially closed or simply $n$-e.c.\ if for any set of vertices $S$ of size $n$ and any set $T\subseteq S$, there is a vertex $x\not\in S$ adjacent to each vertex of $T$ and no vertex of $S\setminus T$. 
We extend this concept to uniform hypergraphs, find necessary conditions for $n$-e.c.\ hypergraphs to exist, and prove that random uniform hypergraphs are asymptotically $n$-existentially closed.
We then provide constructions to generate infinitely many examples of $n$-e.c.\ hypergraphs.
In particular, these constructions use certain combinatorial designs as ingredients, adding to the ever-growing list of applications of designs.
\end{abstract}


MSC Classification Codes: 05B05, 05C65, 05C80, 05C99

\pagebreak
\section{Introduction}

For a positive integer $n$, a graph $G$ with at least $n$ vertices is {\it $n$-existentially closed} or simply {\it $n$-e.c.\ }if for any set of vertices $S$ of size $n$ and any set $T\subseteq S$, there is a vertex $x\in V(G)\setminus S$ adjacent to each vertex of $T$ and no vertex of $S\setminus T$. 
We say that $x$ is {\it correctly joined} to $T$ and $S\setminus T$.
Hence, for each $n$-subset $S$ of vertices, there exist $2^n$ vertices joined to $S$ in all possible ways. 
For example, a 1-e.c.\ graph is one with neither isolated nor universal vertices.

If a graph has the $n$-e.c.\ property, then it possesses other structural properties such as the following.

\begin{Theorem}\cite{Bon}\label{Bon}
Let $G$ be an $n$-e.c.\ graph where $n$ is a positive integer.
\begin{enumerate}
\item The graph $G$ is $m$-e.c.\ for all $1\leq m\leq n-1$.
\item The graph $G$ has order at least $n+2^n$, and has at least $n2^{n-1}$ edges.
\item The complement of $G$ is $n$-e.c.
\item Each graph of order at most $n+1$ embeds in $G$.
\item If $n>1$, then for each vertex $x$ of $G$, each of the graphs $G-x$, the subgraph induced by the neighbourhood $N(x)$, and the subgraph induced by $(V(G)\setminus N(x))-x$ are $(n-1)$-e.c.
\end{enumerate}
\end{Theorem}

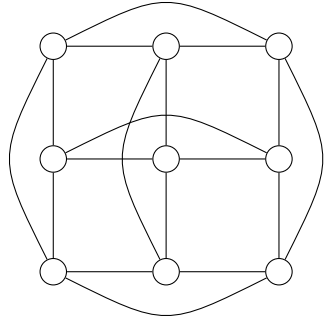
\begin{figure}[h!]
\centering
\vfill
\begin{tikzpicture}[scale=1.5]
\tikzstyle{every node}=[circle, draw, fill=white!75,inner sep=0pt, minimum width=10pt];
\node (v1) at (0,0)[] {};
\node (v2) at (0,1)[] {};
\node (v3) at (0,2)[] {};
\node (v4) at (1,0)[] {};
\node (v5) at (1,1)[] {};
\node (v6) at (1,2)[] {};
\node (v7) at (2,0)[] {};
\node (v8) at (2,1)[] {};
\node (v9) at (2,2)[] {}; 

\draw (v1)--(v2);
\draw (v1).. controls (-0.5,1) ..(v3);
\draw (v2)--(v3);

\draw (v4)--(v5);
\draw (v4).. controls (0.5,1) ..(v6);
\draw (v5)--(v6);

\draw (v7)--(v8);
\draw (v7).. controls (2.5,1) ..(v9);
\draw (v8)--(v9);

\draw (v1)--(v4);
\draw (v1).. controls (1,-0.5) ..(v7);
\draw (v4)--(v7);

\draw (v2)--(v5);
\draw (v2).. controls (1,1.5) ..(v8);
\draw (v5)--(v8);

\draw (v3)--(v6);
\draw (v3).. controls (1,2.5) ..(v9);
\draw (v6)--(v9);

\end{tikzpicture}
\caption{The graph $K_3\square K_3$}\label{K3prod}
\vfill
\end{figure}

Some examples of $n$-e.c.\ graphs include the three non-isomorphic 1-e.c.\ graphs of minimum order 4, depicted in \Cref{2K2,P4,C4}, and the 2-e.c.\ graph of order 9 depicted in Figure~\ref{K3prod}. In~\cite{CEV} it was shown that the minimum order of a 2-e.c.\ graph is nine and in~\cite{BonCam} it was established that $K_3\square K_3$ is the unique 2-e.c.\ graph on nine vertices.

\begin{figure}[h!]
\centering
\begin{minipage}[b][5cm][s]{.31\textwidth}
\centering
\vfill
\begin{tikzpicture}
\tikzstyle{every node}=[circle, draw, fill=white!75,inner sep=0pt, minimum width=7pt];
\node (v1) at (0,0)[] {};
\node (v2) at (0,2)[] {};
\node (v3) at (2,0)[] {};
\node (v4) at (2,2)[] {};

\draw (v1)--(v2);
\draw (v3)--(v4);

\end{tikzpicture}
\caption{The graph $2K_2$}\label{2K2}
\vfill
\end{minipage}\quad
\begin{minipage}[b][5cm][s]{.31\textwidth}
\centering
\vfill
\begin{tikzpicture}
\tikzstyle{every node}=[circle, draw, fill=white!75,inner sep=0pt, minimum width=7pt];
\node (v1) at (0,0)[] {};
\node (v2) at (0,2)[] {};
\node (v3) at (2,0)[] {};
\node (v4) at (2,2)[] {};

\draw (v1)--(v2);
\draw (v1)--(v3);
\draw (v3)--(v4);

\end{tikzpicture}
\caption{The graph $P_4$}\label{P4}
\vfill
\end{minipage}\quad
\begin{minipage}[b][5cm][s]{.31\textwidth}
\centering
\vfill
\begin{tikzpicture}
\tikzstyle{every node}=[circle, draw, fill=white!75,inner sep=0pt, minimum width=7pt];
\node (v1) at (0,0)[] {};
\node (v2) at (0,2)[] {};
\node (v3) at (2,0)[] {};
\node (v4) at (2,2)[] {};

\draw (v1)--(v2);
\draw (v1)--(v3);
\draw (v2)--(v4);
\draw (v3)--(v4);

\end{tikzpicture}
\caption{The graph $C_4$}\label{C4}
\vfill
\end{minipage}
\end{figure}

A {\it hypergraph} $H$ is a pair $(V, E)$ such that $V$ is a set of distinct elements called {\it vertices} and $E$ is a collection of subsets of $V$ called {\it hyperedges} or simply {\it edges}. When every edge in a hypergraph $H$ is of the same cardinality $h$, we say that $H$ is an {\it $h$-uniform} hypergraph.

We extend the notion of an $n$-existentially closed graph to uniform hypergraphs as follows.
For an $h$-uniform hypergraph $H$, we say that $H$ is {\it $n$-e.c.\ }if, for any set of vertices $S$ of size $n$ and any set $T\subseteq S$, there is a set of vertices $X\subseteq V(H)\setminus S$ of size $h-1$ such that for all $z\in T$, $X\cup \{z\}$ is an edge of $H$ and for all $s\in S\setminus T$, $X\cup\{s\}$ is not an edge of $H$.
We again say that the set $X$ is {\it correctly joined} to $T$ and $S\setminus T$.
Note that for $h=2$, this definition agrees with the usual notion of an existentially closed graph.

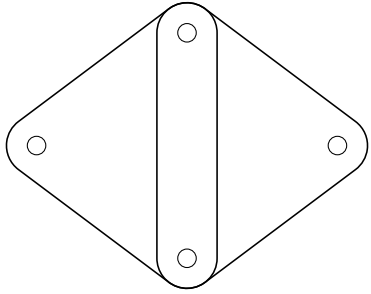
\begin{figure}[h]
\centering
\begin{tikzpicture}
[
    he/.style={draw, semithick},        
    ce/.style={draw, dashed, semithick}, 
]
\tikzstyle{every node}=[circle, draw, fill=white!75,inner sep=0pt, minimum width=7pt];

\node (a) at (0,1.5)[] {};
\node (b) at (2,0)[] {};
\node (c) at (0,-1.5)[] {};
\node (d) at (-2,0)[] {};

\draw[he] \hedgeiii{a}{b}{c}{4mm};
\draw[he] \hedgeiii{d}{a}{c}{4mm};

\end{tikzpicture}\caption{The smallest 1-e.c.\ 3-uniform hypergraph}\label{smallhec}
\end{figure}

For an $h$-uniform hypergraph $H$ to be 1-e.c., for each vertex $x$, the hypergraph must have at least one edge containing $x$ and there must exist at least one set of vertices of size $h-1$ which does not form an edge with $x$.
The smallest example of a 1-e.c.\ $h$-uniform hypergraph can be observed by considering the hypergraph on $h+1$ vertices, with two edges that share exactly $h-1$ vertices in common. For example, Figure~\ref{smallhec} depicts the smallest 1-e.c.\ 3-uniform hypergraph both in the number of edges and vertices.

In Section~\ref{hyperec}, we identify multiple necessary conditions for the existence of existentially closed graphs which extend naturally to existentially closed hypergraphs.
Many of these results mirror those listed in Theorem~\ref{Bon}.
We also prove that random uniform hypergraphs are asymptotically existentially closed.
In particular, for a large enough number of vertices and any $n\geq 1$, random uniform hypergraphs are $n$-existentially closed.
In a sense, this implies that most uniform hypergraphs are existentially closed.
However, as this result is non-constructive, we are still left without examples of such hypergraphs.
In Section~\ref{designs}, we address this situation by presenting constructions for building existentially closed hypergraphs from combinatorial designs.
In particular, we construct infinitely many $n$-e.c.\ uniform hypergraphs for any $n\geq 1$ given appropriate combinatorial designs, which are known to exist whenever the obvious necessary conditions are met and the order is sufficiently large.

\section{$n$-E.C.\ Uniform Hypergraphs}\label{hyperec}

In Section~\ref{necessary}, we present multiple necessary conditions for the existence of $n$-e.c.\ uniform hypergraphs and in Section~\ref{randomsec} we prove that random uniform hypergraphs are asymptotically existentially closed.

\subsection{Necessary Conditions}\label{necessary}

As is the case for graphs, some immediate structural properties of $n$-e.c.\ hypergraphs are easily observed.

\begin{Theorem}\label{smaller}
Let $H$ be an $h$-uniform hypergraph. If $H$ is $n$-e.c., then $H$ is $m$-e.c.\ for each $1\leq m\leq n$.
\end{Theorem}

\begin{Proof}
Let $S=\{v_1,v_2,\dots, v_m\}$ be a set of vertices of $H$ where $1\leq m\leq n$ and choose $T\subseteq S$.
Pick vertices $v_{m+1},v_{m+2},\dots, v_n \in V(H)\setminus S$.
Since $H$ is $n$-e.c.\ there exists an $(h-1)$-set $X$ correctly joined to $T$ and $(S\setminus T) \cup \{v_{m+1},v_{m+2},\dots, v_n\}.$
In particular, this set $X$ is also correctly joined to $T$ and $S\setminus T$ and so $H$ is $m$-e.c.
\end{Proof}

We also identify some lower bounds on the number of vertices and edges in an $n$-e.c.\ hypergraph.

\begin{Theorem}\label{size}
Let $H$ be an $h$-uniform hypergraph. If $H$ is $n$-e.c.\ then $H$ has at least $n2^{n-1}$ edges and at least $n+\ell$ vertices, where $\ell$ is the smallest positive integer such that $\binom{\ell}{h-1}\geq 2^n$.
\end{Theorem}

\begin{Proof}
Let $S$ be an $n$-set in $V(H)$. 
For each $x\in S$, $x$ is contained in $2^{n-1}$ sets $T\subseteq S$, each of which is correctly joined to at least one appropriate set $X$.
Each set $X$ forms an edge with vertex $x$ and thus $\deg(x)\geq 2^{n-1}$.
Thus, since each set $X$ is disjoint from $S$, $H$ has at least $n2^{n-1}$ edges.

Note that $H$ must have at least $n$ vertices, plus enough other vertices to form at least $2^n$ sets of size $h-1$.
Let $\ell$ be the smallest positive integer such that $\binom{\ell}{h-1}\geq 2^n$.
Then $H$ has at least $n+\ell$ vertices.
\end{Proof}
Observe that when $h=2$, $\ell=2^n$ which agrees with part 2 of Theorem 1.

For an $h$-uniform hypergraph $H$, define $H^\mathsf{c}$ as the hypergraph on the vertex set $V(H)$ where an $h$-set $e$ of vertices is an edge of $H^\mathsf{c}$ if and only if $e$ is not an edge of $H$. We call $H^\mathsf{c}$ the {\it $h$-uniform complement} of $H$ or simply, the {\it complement} of $H$.

\begin{Theorem}\label{comp}
Let $H$ be an $h$-uniform hypergraph. If $H$ is $n$-e.c.\ then the complement $H^\mathsf{c}$ is also $n$-e.c.
\end{Theorem}

\begin{Proof}
Let $S$ be an $n$-set of vertices in $V(H^\mathsf{c})$ and $T\subseteq S$.
Since $H$ is $n$-e.c.\ there is an $(h-1)$-set $X$ which is correctly joined to $S\setminus T$ and $T$, meaning $X$ forms an edge in $H$ with each vertex of $S\setminus T$ and with no vertex of $T$.
But this means that $X$ forms an edge in $H^\mathsf{c}$ with each vertex of $T$ and with no vertex of $S\setminus T$.
Thus, $H^\mathsf{c}$ is $n$-e.c.
\end{Proof}

For a hypergraph $H=(V,E)$, we say that a hypergraph $H'=(V',E')$ is a {\it subgraph} of $H$ if $V'\subseteq V$ and $E'\subseteq E$.
In other words, $H'$ is a subgraph of $H$ if and only if every vertex of $H'$ is also a vertex of $H$ and every edge of $H'$ is also an edge of $H$.
Note that this definition of a subgraph of a hypergraph coincides with the definition of a ``strong subhypergraph'' in~\cite{DPP} and a ``hypersubgraph'' in~\cite{BS}.

Now let $H$ be a hypergraph and let $Y\subseteq V(H)$ be a subset of vertices of $H$.
We denote the {\it subgraph induced by $Y$} in $H$ by $H[Y]$.
That is, $H[Y]$ is the hypergraph on the vertex set $Y$ whose edges are precisely the edges of $H$ in which each vertex is a member of $Y$.
Also, for $v\in V(H)$ the {\it neighbourhood} of $v$, denoted $N(v)$, is the set of all vertices which occur together with $v$ in at least one edge of $H$.

\begin{Theorem}\label{induced}
Let $H$ be an $h$-uniform hypergraph. If $H$ is $n$-e.c.\ then for each vertex $v\in~V(H)$, the hypergraphs $H-v$ and $H[N(v)]$ are $(n-1)$-e.c.
\end{Theorem}

\begin{Proof}
Let $S$ be an $(n-1)$-set of vertices in $V(H-v)$ and $T\subseteq S$.
Since $H$ is $n$-e.c.\ there is an $(h-1)$-set $X$ which is correctly joined to $T$ and $(S\cup \{v\})\setminus T$ in $H$.
Note that by definition, $v\not\in X$.
So $X$ forms an edge in $H-v$ with each vertex of $T$ and with no vertex of $S\setminus T$.
Thus, $H-v$ is $(n-1)$-e.c.

Now let $S$ be an $(n-1)$-set of vertices in $H[N(v)]$ and $T\subseteq S$.
Since $H$ is $n$-e.c.\ there is an $(h-1)$-set $X$ which is correctly joined to $T\cup \{v\}$ and $(S\cup \{v\})\setminus (T\cup \{v\})$ in $H$.
Note that since $X$ forms an edge with each vertex of $T\cup \{v\}$, then $X\cup \{v\}$ is an edge of $H$ and so each vertex of $X$ is contained in $N(v)$.
Therefore $X$ forms an edge in $H[N(v)]$ with each vertex of $T$ and with no vertex of $S\setminus T$.
Thus, $H[N(v)]$ is $(n-1)$-e.c.
\end{Proof}

In Section~\ref{designs} we will see examples of $h$-uniform $n$-e.c.\ hypergraphs $H$ constructed from combinatorial designs.
These hypergraphs have the additional property that each pair of vertices appears together in at least one edge of $H$.
This means that for any $v\in V(H)$, the neighbourhood of $v$ is $N(v)=V(H)\setminus \{v\}$ and so the set of non-neighbours of $v$ is $(V(H)\setminus N(v))\setminus \{v\}$ which is an empty set of vertices.
So the subgraph induced by this set is the empty graph and is therefore not $n$-e.c.\ for any $n$.

To continue extending results listed in Theorem \ref{Bon} to $h$-uniform hypergraphs, we define a slightly altered notion of a set of non-neighbours of a vertex in a hypergraph.
For $v\in V(H)$, let $A(v)$ be the set of all vertices that occur together with $v$ in at least one edge of $H^\mathsf{c}$.
Note that for a graph, $A(v)=(V(H)\setminus N(v))\setminus \{v\}$.
With this distinction, we may establish the following result.

\begin{Theorem}\label{Av}
Let $H$ be an $h$-uniform hypergraph. If $H$ is $n$-e.c.\ then for each vertex $v\in~V(H)$, the hypergraph $H[A(v)]$ is $(n-1)$-e.c.
\end{Theorem}

\begin{Proof}
Let $S$ be an $(n-1)$-set of vertices in $H[A(v)]$ and $T\subseteq S$.
Since $H$ is $n$-e.c.\ there is an $(h-1)$-set $X$ which is correctly joined to $T$ and $(S\cup \{v\})\setminus T$ in $H$.
Note that $X$ does not form an edge with any vertex of $(S\cup \{v\})\setminus T$.
In particular, $X$ does not form an edge with $v$, so $X\subseteq A(v)$.
Now note that $X$ forms an edge in $H[A(v)]$ with each vertex of $T$ and with no vertex of $S\setminus T$.
Thus, $H[A(v)]$ is $(n-1)$-e.c.
\end{Proof}

In the next section, we prove that random uniform hypergraphs are asymptotically existentially closed, a result that mirrors that of existentially closed graphs.

\subsection{Random $n$-E.C.\ Hypergraphs}\label{randomsec}

One of the earliest results on existentially closed graphs is that random finite graphs are asymptotically $n$-e.c.\ \cite{ER2}.
We show that this result also extends to $n$-e.c.\ hypergraphs.

A {\it random $h$-uniform hypergraph}, denoted $H_h(m,p)$, is an $h$-uniform hypergraph on $m$ vertices in which each set of vertices $e\subseteq V(H)$ of size $h$ is chosen to be an edge of $H$ randomly and independently with probability $p$, where $p$ may depend on $m$.
Thus, for $h=2$ this model reduces to the well-known Erd\H{o}s-R\'{e}nyi model $G(m,p)$ \cite{ER}.
For some early results on random hypergraphs, see~\cite{KL}.

\begin{Theorem}\label{random}
Let $p$ be a fixed real number such that $0<p<1$ and let $n>1$ and $h>1$ be integers. With probability 1 as $m\rightarrow \infty$, $H_h(m,p)$ satisfies the $n$-e.c.\ property. 
\end{Theorem}

\begin{Proof}
Fix an $n$-set $S$ of vertices and fix $T\subseteq S$. 
For a given $(h-1)$-set $X\subseteq V\setminus S$, the probability that $X$ is not correctly joined to $T$ and $S\setminus T$ is $1-p^n$.
The probability that no set of size $h-1$ is correctly joined to $T$ and $S\setminus T$ is therefore $$\left(1-p^n\right)^{\binom{m-n}{h-1}}.$$ 
As there are $\binom{m}{n}$ choices for $S$ and $2^n$ choices for $T\subseteq S$, the probability that $H_h(m,p)$ is not $n$-e.c.\ is at most $$\binom{m}{n}2^n\left(1-p^n\right)^{\binom{m-n}{h-1}}.$$
Since $n$, $p$, and $h$ are fixed, the probability that $H_h(m,p)$ is not $n$-e.c.\ tends to 0 as $m\rightarrow \infty$.
\end{Proof}

Since random uniform hypergraphs asymptotically satisfy the $n$-e.c.\ property, we should expect to see many examples of $n$-e.c.\ uniform hypergraphs.
However, as is the case for graphs, it is not immediately clear how to find examples of these hypergraphs.
In the next section, we detail constructions for building existentially closed hypergraphs from combinatorial designs, namely, Latin squares, balanced incomplete block designs, and $t$-designs.

\section{$n$-E.C.\ Hypergraphs from Designs}\label{designs}

Since we know that random uniform hypergraphs are asymptotically $n$-e.c.\ for any $n\geq 1$, we now look for constructions that generate infinite families of such hypergraphs.
One such construction for generating 2-e.c.\ uniform hypergraphs makes use of a set of well-known objects within combinatorics, Latin squares.

A {\it Latin square of order $n$} is an $n\times n$ array consisting of $n$ symbols such that each symbol occurs in each row and each column precisely once.
A pair of Latin squares of the same order is said to be {\it orthogonal} if, when superimposed, the entries viewed as ordered pairs are all unique.
A set of Latin squares of the same order in which any two form an orthogonal pair is said to be a set of {\it mutually orthogonal Latin squares} or {\it MOLS} for short.
It is well-known that the maximum possible number of mutually orthogonal Latin squares of order $n$ is $n-1$.
Such a set of MOLS is referred to as a {\it complete set} of MOLS.
Complete sets of MOLS of order $n$ are known to exist when $n$ is a prime or power of a prime.
For more information on Latin squares, including applications, see~\cite{DK}.
For an example of a complete set of MOLS of order 4, see Figure~\ref{MOLS4}.

\vspace{12pt}
\begin{figure}[h!]
\centering
\begin{tabular}{|c|c|c|c|}
\hline
0&1&2&3\\
\hline
2&3&0&1\\
\hline
1&0&3&2\\
\hline
3&2&1&0\\
\hline
\end{tabular}
\qquad
\begin{tabular}{|c|c|c|c|}
\hline
0&1&2&3\\
\hline
3&2&1&0\\
\hline
2&3&0&1\\
\hline
1&0&3&2\\
\hline
\end{tabular}
\qquad
\begin{tabular}{|c|c|c|c|}
\hline
0&1&2&3\\
\hline
1&0&3&2\\
\hline
3&2&1&0\\
\hline
2&3&0&1\\
\hline
\end{tabular}
\caption{A complete set of MOLS of order 4}\label{MOLS4}
\end{figure}

Now suppose $L$ is a set of $\ell$ mutually orthogonal Latin squares of order $h+1$. 
We will form an $h$-uniform hypergraph $H_L$ in the following way.
Let $V$ be a set of $(h+1)^2$ vertices organised into a $(h+1)\times(h+1)$ array $A$.
We define edges to form the edge set $E$ in two ways.
Firstly, for each row (respectively each column) of $A$, take all $(h+1)$ $h$-sets of vertices within the row (respectively column) to be edges of $E$.
Secondly, for each Latin square in $L$ and for each symbol within the squares of $L$, take note of the position of each occurrence of that symbol and then take the corresponding $(h+1)$-set of vertices within $A$. 
Now take all $h$-subsets of vertices within this set as edges of $E$.
The resulting hypergraph $H_L=(V,E)$ is an $h$-uniform hypergraph on $(h+1)^2$ vertices and $(\ell+2)(h+1)^2$ edges.

\vspace{12pt}
\begin{figure}[h!]
\centering
\begin{tabular}{|c|c|c|c|}
\hline
0&1&2&3\\
\hline
4&5&6&7\\
\hline
8&9&a&b\\
\hline
c&d&e&f\\
\hline
\end{tabular}
\caption{The $4\times 4$ array $A$}\label{array}
\end{figure}

For example, consider the complete set of MOLS in Figure~\ref{MOLS4} as our ingredient set $L$. 
Then the constructed hypergraph $H_L=(V,E)$ is a 3-uniform hypergraph on 16 vertices.
Organise the vertex set $V$ into a $4\times 4$ array $A$ (see Figure~\ref{array}).
Then form the edge set $E$ according to the construction above. 
For instance, the 4-sets acquired from the symbols in the first square in $L$ are the following: $$\{0,6,9,f\}, \{1,7,8,e\}, \{2,4,b,d\}, \{3,5,a,c\}.$$
We then take all 3-subsets of these 4-sets as edges of $E$. 
The resulting hypergraph $H_L$ will have 80 edges according to our construction. 

Our next theorem asserts that if $L$ is a complete set of MOLS, then the resulting hypergraph $H_L$ is 2-e.c.

\begin{Theorem}\label{MOLS}
If $L$ is a complete set of MOLS of order $h+1$ and $h\geq 3$, then the hypergraph $H_L$ is 2-e.c.
\end{Theorem}

\begin{Proof}
We must verify that for any 2-set of vertices $S=\{u,v\}$ and any $T\subseteq S$, there is an $(h-1)$-set $X\subseteq V(H_L)\setminus S$ such that $X$ forms an edge with each vertex of $T$ and with no vertex of $S\setminus T$.

When $|T|=0$, we take $X$ to be any $(h-1)$-set of vertices all of which occur in the same row of $A$ together but a different row than that of $u$ and $v$. 
Here, $X$ forms an edge with neither $u$ nor $v$.
When $|T|=1$, say $T=\{u\}$, if $u$ and $v$ are in distinct columns of $A$ then we take $X$ to be an $(h-1)$-set of vertices which each occur in the same column as $u$ and note that $X$ forms an edge with $u$ but not with $v$.
Otherwise, if $u$ and $v$ are in the same column, we take $X$ to be an $(h-1)$-set of vertices which each occur in the same row of $A$ as $u$ and note that $X$ forms an edge with $u$ but not with $v$.

Finally, when $|T|=2$, if $u$ and $v$ happen to be in the same row (respectively column) of $A$, we take $X$ to be the set containing the other $h-1$ vertices in that row (respectively column) and note that $X$ forms an edge with each of $u$ and $v$.
Otherwise, when $u$ and $v$ are in distinct rows and columns, we find the Latin square among $L$ in which the positions corresponding to $u$ and $v$ share a common symbol; this is guaranteed to exist since $L$ is a complete set of MOLS.
We then take $X$ to be the set containing the other $h-1$ vertices corresponding to the positions of the other occurrences of the shared symbol in that Latin square and note that $X$ forms an edge with each of $u$ and $v$.
\end{Proof}

Since complete sets of MOLS are known to exist whenever the order of the Latin squares is a prime or prime power, Theorem~\ref{MOLS} implies that there are infinitely many 2-e.c.\ uniform hypergraphs.

We now generalise our construction to produce an infinite family of 2-e.c.\ uniform hypergraphs, this time using a wider set of combinatorial objects as our initial ingredients, balanced incomplete block designs.

A {\it balanced incomplete block design} or {\it BIBD} with parameters $(v,k,\lambda)$ is a pair $\mathcal{D}=(V,\mathcal{B})$ such that $V$ is a set of $v$ distinct elements called {\it points} and $\mathcal{B}$ is a collection of $k$-subsets of $V$ called {\it blocks} such that each pair of points of $V$ occurs in exactly $\lambda$ blocks of $\mathcal{B}$.
The total number of blocks in a BIBD is denoted $b$ and the number of blocks which contain any given point is called the {\it replication number} and is denoted $r$.
In particular, $b= \frac{\lambda v(v-1)}{k(k-1)}$ and $r= \frac{\lambda(v-1)}{k-1}$.
Balanced incomplete block designs are known to exist asymptotically whenever the necessary divisibility conditions are met~\cite{W1,W2,W3}.
For more information on BIBDs, including constructions and known examples, see~\cite{HCD}.

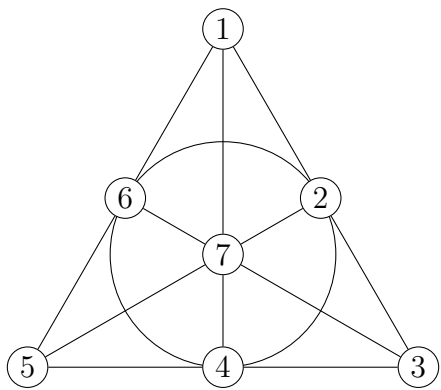
\begin{figure}[h]
    \centering
    \begin{tikzpicture}[scale=1.5]
\tikzstyle{every node}=[circle, draw, fill=white!75,inner sep=2pt, minimum width=7pt];
\node (v7) at (0,0) [] {7};
\draw (0,0) circle (1cm);
\node (v1) at (90:2cm) [] {1};
\node (v2) at (210:2cm) [] {5};
\node (v4) at (330:2cm) [] {3};
\node (v3) at (150:1cm) [] {6};
\node (v6) at (270:1cm) [] {4};
\node (v5) at (30:1cm) [] {2};
\draw (v1) -- (v3) -- (v2);
\draw (v2) -- (v6) -- (v4);
\draw (v4) -- (v5) -- (v1);
\draw (v3) -- (v7) -- (v4);
\draw (v5) -- (v7) -- (v2);
\draw (v6) -- (v7) -- (v1);
\end{tikzpicture}
    \caption{The Fano plane}
    \label{fano}
\end{figure}

For example, Figure~\ref{fano} is a graphical representation of the unique (7,3,1)-BIBD.
This design is known as the Fano plane and is significant within multiple branches of mathematics including design theory, projective geometry, and group theory. 
Here, each line (including the curved one) passes through exactly three points and determines a block of size 3 consisting of those points. 
The corresponding blocks of this design are then $\{1,2,3\}, \{3,4,5\},\{1,5,6\},\{1,4,7\},\{2,5,7\},\{3,6,7\},
\text{ and }\{2,4,6\}.$
Note that each pair of points occurs in precisely one block.
For more information on the Fano plane see~\cite{HCD}.

Any design with $v>k$, other than the design in which the block set $\mathcal{B}$ is the set of all $k$-subsets of $V$, is actually 1-e.c.\ when viewed as a $k$-uniform hypergraph.
Indeed, each point occurs in exactly $r$ blocks, so each point occurs at least once and no point occurs $\binom{v}{k}$ times.

Now suppose $\mathcal{D}$ is a $(v,k,1)$-BIBD with $k\geq 3$. 
For each $h$ such that  $3\leq h\leq k$, we will form an $h$-uniform hypergraph $H_{\mathcal{D},h}$ in the following way.
Let the vertex set $V$ be the point set of $\mathcal{D}$.
For each block $B$ of $\mathcal{D}$, take all $\binom{k}{h}$ $h$-subsets of $B$ as edges of the edge set $E$.
The resulting hypergraph $H_{\mathcal{D},h}=(V,E)$ is an $h$-uniform hypergraph with $v$ vertices and $b \binom{k}{h}$ edges where $b$ is the total number of blocks in $\mathcal{D}$.

\begin{Theorem}\label{BIBD}
Let $\mathcal{D}$ be a $(v,k,1)$-BIBD with $k\geq 4$.
If $v\geq k+2$ and $3\leq h\leq k-1$, then the hypergraph $H_{\mathcal{D},h}$ is 2-existentially closed.
\end{Theorem}

\begin{Proof}
We must verify that for any 2-set of vertices $S=\{u,v\}$ and any $T\subseteq S$, there is an $(h-1)$-set $X\subseteq V(H_{\mathcal{D},h})\setminus S$ such that $X$ forms an edge with each vertex of $T$ and with no vertex of $S\setminus T$.
Note that since $h\geq 3$ and $\lambda=1$, the number of times any $(h-1)$-set $X$ occurs within a block of $\mathcal{D}$ is at most once; otherwise, the block containing $X$ would contain a pair which occurs more than $\lambda$ times among the blocks of $\mathcal{D}$.
So any $(h-1)$-set $X$ chosen directly from a block of $\mathcal{D}$ is unique.

When $|T|=2$, let $B$ be the unique block containing both $u$ and $v$. 
Since $k\geq h+1$, $B$ contains $u$, $v$ and at least $h-1$ other points.
So take $X$ to be an $(h-1)$-set of points in $B\setminus \{u,v\}$ and note that $X$ forms an edge with each of $u$ and $v$ in $H_{\mathcal{D},h}$.

When $|T|=1$, say $T=\{u\}$, let $B$ be a block among the $r$ blocks that contain $u$ other than the unique block which contains both $u$ and $v$.
Such a block exists since the replication number $r= \frac{v-1}{k-1}$ is greater than 1 whenever $v>k$.
Now take $X$ to be an $(h-1)$-set of points other than $u$ in this block and note that $X$ forms an edge with $u$ but not with $v$.

Finally, when $|T|=0$, if there exists a block $B$ which contains neither $u$ nor $v$, then we can choose an $(h-1)$-set $X$ within $B$ and note that $X$ forms an edge with neither $u$ nor $v$.
There are $b$ blocks in total, $r$ blocks containing $u$, $r$ blocks containing $v$, and exactly one block containing both $u$ and $v$ (which is counted twice among the blocks containing $u$ and $v$).
So if $b>2r-1$ then there exists an appropriate block $B$ from which to choose an $(h-1)$-set $X$.
Recall that $b= \frac{\lambda v(v-1)}{k(k-1)}$, $r= \frac{\lambda(v-1)}{k-1}$ and $\lambda=1$, so
$$\begin{array}{lrcl}
 &\frac{v(v-1)}{k(k-1)} &>& \frac{2(v-1)}{k-1} -1    \\
\Leftrightarrow &v(v-1) &>& 2k(v-1) - k(k-1)      \\
\Leftrightarrow &(v-2k)(v-1) &>& k(1-k).
\end{array}$$
Now when $v\geq k+2$, $$(v-2k)(v-1)\geq (2-k)(k+1) > k(1-k)$$ and thus $b>2r-1$ holds.
\end{Proof}

Note that if $h=2$ this construction would yield a graph (i.e., a 2-uniform hypergraph) but such a graph would not even be 1-existentially closed.
Indeed, since each pair of points in a design occurs precisely $\lambda$ times, the resulting graph would be complete and trivially not 1-e.c.
Also, if $h=k$ then the hypergraph $H_{\mathcal{D},h}$ is simply the design $\mathcal{D}$ itself.
Note that any design with $\lambda=1$ when viewed as a hypergraph cannot be 2-e.c.\ since by definition, there would need to exist a set $X$ of size $k-1$ which forms an edge (or block) with at least two distinct points, violating $\lambda=1$.
However, finding designs with higher values of $\lambda$ that are $n$-e.c.\ for $n\geq 2$ is an open problem.

To find examples of $n$-e.c.\ hypergraphs for values of $n\geq 3$, we make use of a natural generalisation of balanced incomplete block designs.
A {\it $t$-$(v,k,\lambda)$ block design}, or {\it $t$-design} for short, is a pair $\mathcal{D}=(V,\mathcal{B})$ such that $V$ is a set of $v$ distinct points and $\mathcal{B}$ is a collection of blocks of size $k$ such that each $t$-subset of points of $V$ occurs in exactly $\lambda$ blocks of $\mathcal{B}$.
Note that a $t$-design with $t=2$ is exactly a balanced incomplete block design.
Infinitely many nontrivial $t$-designs without repeated blocks are known to exist for all $t$ \cite{Teir}.
Asymptotically, $t$-designs are known to exist whenever the necessary divisibility conditions are met~\cite{Keevash}.
For more information on $t$-designs including constructions and known examples, see~\cite{HCD}.

Now suppose $\mathcal{D}$ is a $t$-$(v,k,1)$-design with $k\geq 3$. 
For each $h$ such that $3\leq h\leq k$, we will form an $h$-uniform hypergraph  $H_{\mathcal{D},h}$ in the following way.
Let the vertex set $V$ be the point set of $\mathcal{D}$.
For each block $B$ of $\mathcal{D}$, take all $\binom{k}{h}$ $h$-subsets of $B$ as edges of the edge set $E$.
The resulting hypergraph  $H_{\mathcal{D},h}=(V,E)$ is an $h$-uniform hypergraph with $v$ vertices and $b \binom{k}{h}$ edges where $b$ is the number of blocks in $\mathcal{D}$.

To show that $H_{\mathcal{D},h}$ is an existentially closed hypergraph for certain values of $v$, $k$, and $h$, we make use of a result that can be found as a remark in Chapter II, Section 4.2 of~\cite{HCD} that allows us to count the number of blocks in a design that contain certain points while avoiding other certain points.

\begin{Lemma}\cite[\S II.4.2]{HCD}\label{lemma}
Let $\mathcal{D}=(V,\mathcal{B})$ be a $t$-$(v,k,\lambda)$ block design and let $I$ and $J$ be disjoint subsets of $V$ with $|I|=i$, $|J|=j$ and $i+j\leq t$.
If $\lambda_{i,j}$ is the number of blocks that contain each point of $I$ and no point of $J$, then $\lambda_{i,j}=\lambda\binom{v-i-j}{k-i}/\binom{v-t}{k-t}$.
\end{Lemma}

\begin{Theorem}\label{t-design}
Let $\mathcal{D}$ be a $t$-$(v,k,1)$-design with $k\geq 2t$. 
If $v\geq k+t$ and $t+1\leq h\leq k-t+1$, then the hypergraph $H_{\mathcal{D},h}$ is $t$-existentially closed.
\end{Theorem}

\begin{Proof}
We must verify that for any $t$-set of vertices $S$ and any $T\subseteq S$, there is an $(h-1)$-set $X\subseteq V(H_{\mathcal{D},h})\setminus S$ such that $X$ forms an edge with each vertex of $T$ and with no vertex of $S\setminus T$.
Note that since $h\geq t+1$ and $\lambda=1$, the number of times any $(h-1)$-set $X$ occurs within a block of $\mathcal{D}$ is at most once; otherwise, the block containing $X$ would contain a $t$-set which occurs more than $\lambda$ times among the blocks of $\mathcal{D}$.
So any $(h-1)$-set $X$ chosen directly from a block of $\mathcal{D}$ is unique.

When $|T|=t$, let $B$ be the unique block containing all $t$ points of $T$.
Since $k\geq t+h-1$, $B$ contains the $t$ points of $T$ and at least $h-1$ other points.
So take $X$ to be an $(h-1)$-set of points in $B\setminus T$ and note that $X$ forms an edge with each vertex of $T$ in $H_{\mathcal{D},h}$.

Now suppose $0\leq |T|\leq t-1$.
For notational simplicity, let $|T|=i$.
If there exists a block $B$ which contains the $i$ points of $T$ but none of the $t-i$ points of $S\setminus T$, then we can choose an $(h-1)$-set $X$ consisting of points of $B$ other than the $i$ points of $T$ and note that $X$ forms an edge with each vertex of $T$ but with no vertex of $S\setminus T$.
By Lemma~\ref{lemma}, the number of such blocks is precisely $\lambda_{i,t-i}=\binom{v-t}{k-i}/\binom{v-t}{k-t}$.
Note that $\lambda_{i,t-i}$ is a positive integer so long as $v-t\geq k-i$ and $v-t\geq k-t$.
Since $0\leq i\leq t-1$, these inequalities hold whenever $v\geq k+t$.
\end{Proof}

Since infinitely many $t$-designs are known to exist for all $t$ \cite{Teir}, Theorem~\ref{t-design} implies that there are infinitely many $n$-e.c.\ hypergraphs for any $n$.

\section{Open Problems}

In Section \ref{designs}, we noted that all non-trivial balanced incomplete block designs are 1-e.c.\ when viewed as hypergraphs. 
However, if a BIBD $\mathcal{D}$ has $\lambda=1$, then $\mathcal{D}$ viewed as a hypergraph cannot be 2-e.c.\ since by definition, there would need to exist a set $X$ of size $k-1$ which forms an edge (or block) with at least two distinct points, violating $\lambda=1$.
In general, if a BIBD viewed as a hypergraph is $n$-e.c., then $n\leq\lambda$.
\begin{Question}
Which BIBDs (if any) viewed as hypergraphs are $n$-e.c.\ for $n\geq 2$?
\end{Question}
\noindent We can ask a similar question for $t$-designs and in general, other types of designs such as cycle systems and group divisible designs.
\begin{Question}
Which designs viewed as hypergraphs are existentially closed?
\end{Question}

One explicit family of $n$-e.c.~graphs is the set of Paley graphs.
Paley graphs are constructed from finite fields by taking the set of elements of a field as the vertex set and making two vertices adjacent if those elements differ by a quadratic residue in the field.
For more information on Paley graphs see \cite{GR}.
In~\cite{BEH} and~\cite{BT} it was shown that for any $n$, every sufficiently large Paley graph is $n$-e.c.

One possible area of future work involves investigating whether {\it Paley hypergraphs} satisfy the existential closure property.
Fortunately, the notion of a Paley hypergraph already exists in the literature.
This extension of Paley graphs was first introduced by Kocay in~\cite{kocay} and was later refined by Poto\v{c}nik and \v{S}ajna in~\cite{PoSa} and again by Dueck (Gosselin) in~\cite{dueck}.
\begin{Question}
Under what conditions are Paley hypergraphs existentially closed?
\end{Question}

Additionally, as our definition of existentially closed hypergraphs is specific to uniform hypergraphs, it remains to see if there is a similar concept for non-uniform hypergraphs.
\begin{Question}
What is an appropriate notion of an existentially closed non-uniform hypergraph?
\end{Question}

\section{Acknowledgements}

Authors Burgess and Pike acknowledge NSERC Discovery Grant support and Luther acknowledges NSERC scholarship support.

\end{doublespace}

\end{document}